\documentclass[11pt]{article}

\usepackage{shapepar,bm}
\usepackage[dvips,arrow,matrix,ps,color,line]{xy}
\usepackage{theorem,amsfonts,amssymb,amscd}
\usepackage{geometry}
\usepackage{graphicx}
\usepackage{graphics}
\usepackage{amssymb}
\usepackage{epstopdf}

\usepackage{amssymb, graphics}
\usepackage[colorlinks=true, pdfstartview=FitV, linkcolor=blue, citecolor=blue, urlcolor=blue]{hyperref}

\newcommand{\mathsym}[1]{{}}

\usepackage[usenames]{color}
\definecolor{MyLightMagenta}{cmyk}{0.1,0.8,0,0.1}
\definecolor{MyDarkBlue}{rgb}{0.1,0,0.3}

\makeindex

\def\NN{\mathbb N}

\def\ZZ{\mathbb Z}
\def\CC{\mathbb C}

\def\PP{\mathbb P}

\def\TT{{\bf T}}

\def\SS{{\tilde{S}}}

\def\ttp{{\tt p}}

\def\ikformu{{\mu^{i_1}\wedge\ldots\wedge\mu^{i_k}}}

\def\kformu{{\mu^{1}\wedge\ldots\wedge\mu^{k}}}

\def\ep{{\epsilon}}

\def\ev{{\rm ev}}

\def\Bw{{\bigwedge}}

\def\w2M{\bigwedge^2M}

\def\wM{\bigwedge M}

\def\w{\wedge }
\def\bw{\bigwedge }

\protect\def\wkM{{\bigwedge^kM}} 
\protect 
\protect
\protect

\protect
\protect

\def\sra{\rightarrow}
\def\lra{\longrightarrow}

\def\proof{\noindent{\bf Proof.}\,\,}
\def\qed{{\hfill\vrule height4pt width4pt depth0pt}\medskip}
\def\be{\begin{equation}}
\def\ee{\end{equation}}
\def\bclm{\begin{claim}}
\def\eclm{\end{claim}}
\def\beqn{\begin{eqnarray}}
\def\eeqn{\end{eqnarray}}
\def\beqn*{\begin{eqnarray*}}
\def\eeqn*{\end{eqnarray*}}

\def\kformep{{\epsilon^1\wedge\ldots\wedge\epsilon^k}}

\def\ikformep{\epsilon^{i_1}\wedge\ldots\wedge\epsilon^{i_k}}

\theoremstyle{change}
\theorembodyfont{\rmfamily}

\newtheorem{claim}{}[section]

\global
\global

\def\no@breaks#1{{\def\\{ \ignorespaces}#1}}    

  \makeatletter
\def\cleardoublepage{\clearpage\if@twoside \ifodd\c@page\else
\hbox{} \thispagestyle{empty}
\newpage
\if@twocolumn\hbox{}\newpage\fi\fi\fi} \makeatother

\usepackage{eso-pic,graphicx}
\makeatletter
\newcommand\BackgroundPicture[2]{%
  \setlength{\unitlength}{1pt}%
  default \put(0,\strip@pt\paperheight){%
  \parbox[t][\paperheight]{\paperwidth}{%
    \vfill
     \centering \includegraphics[angle=#2, width=15cm, height=15cm,  bb=0 0 150 150]{#1}
    \vfill
}}} %
\makeatother

\usepackage{geometry}                
\usepackage{graphicx}
\usepackage{amssymb}
\usepackage{epstopdf}

\title{Equivariant Schubert Calculus\thanks{Work partially sponsored by PRIN ``Geometria sulle Variet\`a Algebriche" (Coordinatore Alessandro Verra), INDAM-GNSAGA and Scuola di Dottorato (ScuDo) del Politecnico di Torino.}}
\author{Letterio Gatto, Ta\'\i se Santiago}
\date{}                                           

\begin{document}
\maketitle
\abstract{\small
\noindent
Let $T$ be a torus acting on $\CC^n$ in such a way that, for all $1\leq k\leq n$, the induced action on the grassmannian $G(k,n)$ has only isolated fixed points.  The $T$-equivariant cohomology  of $G(k,n)$ can be realized, in this case, as the quotient of a ring generated by derivations on the exterior algebra of a free module of rank $n$ over the $T$-equivariant cohomology of a point. This result allows  a simple  explicit description of $T$-equivariant Schubert calculus on grassmannians.
}

\section{Introduction}

\claim\label{int01} {\bf Equivariant Schubert Calculus.} Let $T$ be an algebraic or compact torus acting on $\CC^n$ is such a way that the following condition holds:

\smallskip
\noindent
($*$)  {\em for each $1\leq k\leq n$, the induced action on $G(k,n)$, the grassmannian variety parametrizing $k$-dimensional vector subspaces of $\CC^n$,  has only isolated fixed points.}

\smallskip
\noindent
In this case $H^*_T(G(k,n))$,  the $T$-equivariant cohomology ring of $G(k,n)$, is a free $A$-module of finite rank, where $A:=H_T^*(pt)$ is  the cohomology ring of a point. If the torus $T$ has dimension $d$,  then $A$ is  the polynomial ring $\ZZ[y_1,\ldots, y_d]$ in $d$-indeterminates.
 Let  $\{\sigma_I: I\in{\cal I}^k_n\}$ be any  $A$-basis of $H_T^*(G(k,n))$, where ${\cal I}^k_n$ is the set of all increasing sequences of $k$ positive integers not bigger than $n$ (as in~\cite{anderson}, \S 5, Section~1). One can then speak  of ``equivariant Schubert calculus", which concerns the structure constants $C^K_{IJ}\in A$ in the product expansion $\sigma_I\cdot \sigma_J=\sum_{K\in{\cal I}^k_n}C^K_{IJ}\sigma_K$. 

\claim{\bf The result.} \label{seint02}$T$-equivariant Schubert calculus on $\PP^{n-1}:=G(1,n)$ is easy. In fact $H_T^*(\PP^{n-1})$ is a free $A$-module generated by $1,\xi,\ldots, \xi^{n-1}$, where $\xi=c_1^T(O_{\PP^{n-1}}(1))$ is the first equivariant Chern class of the hyperplane bundle over $\PP^{n-1}$. The product structure is then the same as that of the ring $A[X]/\ttp$, where $\ttp$ (following~\cite{LakTh},~\cite{LakTh1}) is the minimal polynomial of the ``multiplying-by-$\xi$" endomorphism of $H_T^*(\PP^{n-1})$. 
Let $(M(\ttp),\jmath)$ be a pair where $M(\ttp)$ is a free $A$-module of rank $n$ and  $\jmath:H_T^*(\PP^{n-1})\sra M(\ttp)$ is an $A$-module isomorphism making $M(\ttp)$ into a free $H_T^*(\PP^{n-1})$-module of rank $1$. Let $(\mu^1,\ldots,\mu^n)$ be any $A$-basis of $M(\ttp)$ such that $\jmath^{-1}(\mu^i)\in H_T^{2i}(\PP^{n-1})$. Then $\mu^1$ generates $M(\ttp)$ as a module over $H_T^*(\PP^{n-1})$. For each $k\geq 1$, we construct an $A$-algebra ${\cal A}^*(\wkM(\ttp))$ such that i) there is an $A$-module isomorphism $\Pi_k:{\cal A}^*(\wkM(\ttp))\sra \wkM(\ttp)$ making $\wkM(\ttp)$ into a free ${\cal A}^*(\wkM(\ttp))$-module generated by $\kformu$; ii) there is a certain sequence $D:=(D_0,D_1,\ldots)$  of derivations on $\wM(\ttp)$ (as in~\ref{rec24}) such that ${\cal A}^*(\wkM(\ttp))$ is a free $A$-module generated by $\{\Delta_I(D)\,|\, I\in {\cal I}^k_n\}$, where, if $I:=(i_1,\ldots,i_k)\in{\cal I}^k_n$,  $\Delta_I(D)$ is the {\em Schur polynomial} $\det((D_{i_j-i})_{1\leq i,j\leq k})$. Our main result, Theorem~\ref{2ndmint}, shows that {\em for each $k\geq 1$ there is a canonical ring isomorphism $\iota_k:H_T^*(G(k,n))\sra {\cal A}^*(\wkM(\ttp))$}. In particular $\jmath_k:=\iota_k\circ\Pi_k$ is an isomorphism, such that $\jmath_1=\jmath$, making $\wkM(\ttp)$ into a free $H_T^*(G(k,n))$-module of rank $1$ generated by $\kformu$. The fact that $H_T^*(G(k,n))$ can be seen as an $A$-algebra (finitely) generated by derivations of $\wM(\ttp)$ allows us to reduce ($T$-equivariant) Schubert calculus on $G(k,n)$  to that, much easier,  on $G(1,n)$ (in the same spirit of~\cite{Gat1}; see also~\cite{Gat2}).

\claim{\bf Related Literature.}\label{int03} Various authors have dealt with $T$-equivariant Schubert calculus in the case when

\smallskip
\noindent
($**$) {\em $T$ is a $n$-dimensional torus (either $(\CC^*)^n$ or $(S^1)^n$) acting diagonally on $\CC^n$. }

\smallskip
\noindent
In an important  paper (\cite{KT}), Knutson and Tao apply the beautiful combinatorics of {\em puzzles}  to a natural basis $\SS$ of $H_T^*(G(k,n))$, explicitly  described  via a list of polynomials satisfying certain GKM conditions (after Goresky, Kottwitz and MacPherson; \cite{GKM}).
 Their main result is that {\em puzzles compute equivariant cohomology} (\cite{KT}; Theorem 2). They also prove a formula ({\em the equivariant Pieri rule}, Proposition 2) governing the product of any element of $\SS$ with $\SS_{div}$, the  generator of $H^2_T(G(k,n))$. 
 A quite different approach is pursued in~\cite{LRS}, whose aim is a description formally closer to  classical Schubert calculus. The authors propose there a very interesting extension of Giambelli's formula but leave open the question about how a full set of Pieri formulas should look like. Knutson and Tao, however, remark that these can be deduced by specializing  those  got by Robinson  for complete flag varieties (\cite{robinson}). 
 
The quantum deformation of $T$-equivariant Schubert calculus in the case ($**$)  is supplied by Mihalcea, in \cite{Mihalcea1} and~\cite{Mihalcea2}. He gets a lot of important results, such as a presentation of $QH_T^*(G(k,n))$, in terms of generators and relations,  as well as a very satisfying (quantum) equivariant Giambelli's type determinantal formula. Classical $T$-equivariant Schubert Calculus is recovered by setting to zero the quantum deformation parameter.  In~\cite{LakTh} and~\cite{LakTh1}, Laksov and Thorup provide a very general description  of  Schubert Calculus  in terms of the action of the ring of symmetric functions on the exterior power of a polynomial ring. Very general Giambelli's and Pieri's formulas are already built in the theory. In~\cite{formalism}, Laksov sets explicitly the connection between  the general Schubert calculus exposed  in~\cite{LakTh} and~\cite{LakTh1}, and the work of Knutson and Tao, and  Mihalcea. 
\claim{\bf An Example.} Our  Theorem~\ref{2ndmint} holds for any action of a $d$-dimensional torus  ($d$ not necessarily equal to $n$) on $\CC^n$,  provided that ($*$) of~\ref{int01} is fulfilled. We devote Section~\ref{KTex}  to show how our methods cope with the situation ($**$) of~\ref{int03}, too. Let $\SS_1,\ldots,\SS_n$ be the basis of $H_T^*(G(1,n))$ used in~\cite{KT}, or~\cite{Mihalcea1} and~\cite{Mihalcea2}, and let $\mu^i=\jmath_1(\SS_i)$ (Cf. Section~\ref{seint02}).  We shall compute very explicit ``equivariant" Pieri's formulas (Theorem~\ref{TeoremaPieriEquivariante}) for the basis  $\{\jmath_k^{-1}(\kformu)\,|\, 1\leq i_1<\ldots<i_k\leq n\}$ of  $H_T^*(G(k,n))$,  fully answering to the wishes of~\cite{LRS}. Such a basis (Cf.~\ref{ex55}) will also allow us to recover  the {\em divisorial Pieri formula} as in~\cite{KT}, Proposition 2. Finally  we  revisit (\ref{KTrevisited}) the example regarding the grassmannian $G(2,4)$ of~\cite{KT}, p.~231, showing that the  GKM-conditions come automatically into play  into  the matrix expression of $\jmath_2^{-1}(\mu^i\w\mu^j)\in H_T^*(G(2,4))$ with respect to  the basis $\{\mu^i\w\mu^j\}$ of $\bw^2M(\ttp)$. See also~\cite{tesi} for more details and examples.

\claim{} {\bf About the proof of~\ref{2ndmint}.}  Let $F_k:=O_{G(k,n)}\otimes\CC^n$. Recall that all the bundles occurring into the universal exact sequence over $G(k,n)$, $0\sra S_k\sra F_k\sra Q_k\sra 0$, are $T$-equivariant. It makes  sense, then, to speak of their $T$-equivariant Chern classes.  We first prove (Proposition~\ref{mainlemma}) that {\em  $H_T^*(G(k,n))$ is a free $A$-module generated by $\{\delta_I:=\Delta_I(c^T(Q_k-F_k))\,\,|\,\, I\in{\cal I}^k_n\}$} where, if $I=(i_1,\ldots,i_k)\in{\cal I}^k_n$, $\Delta_I(c^T(Q_k-F_k))=\det(c_{i_j-i}(Q_k-F_k))$ is the corresponding Schur polynomial in the $T$-equivariant Chern classes of $Q_k-F_k\in K^0(G(k,n))$.  This fact is surely  well known to the experts but  we include a proof because we have not been able to find an explicit reference in the literature. 
Since the equivariant cohomology of $G(k,n)$ can be approximated by the ordinary cohomology of a sequence of certain grassmann bundles $G(k,{\cal E}_m)$ over finite dimensional smooth projective varieties $B_mT$ (Cf. \ref{comdiag32}), the basis we use for $H_T^*(G(k,n))$  can be seen as the ``limit" of the basis for $H^*(G(k,{\cal E}_m))$ described in~\cite{Fu1}, p.~266, and is, therefore, the closest to that used in the ordinary cohomology of Grassmann bundles. With respect to such a canonical basis, in particular, equivariant Schubert calculus comes already equipped with Giambelli's and Pieri's type formulas (\cite{Fu1}, p.~266). Moreover, the coefficients $C^K_{IJ}\in A$ occurring  in the product expansion $\delta_I\delta_J=\sum_{I\in{\cal I}^k_n}C^K_{IJ}\delta_K$   are precisely the classical (integral!) Littlewood-Richardson coefficients if $wt(I)+wt(J)\leq k(n-k)$ (where if $I=(i_1,\ldots, i_k)\in{\cal I}^k_n$, one sets $wt(I)=\sum_{j=1}^k(i_j-j)$). The explicit  module isomorphism $\iota_k:H_T^*(G(k,n))\sra {\cal A}^*(\wkM(\ttp))$ given by $\Delta_I(c^T(Q_k-F_k))\mapsto \Delta_I(D)$ (\ref{2ndmint}), turns out to be a ring isomorphism after invoquing (our rephrasing of) the general Schubert calculus developed by Laksov and Thorup in~\cite{LakTh} and~\cite{LakTh1} (see Section~\ref{preliminaries} for details).

\smallskip
\noindent
 {\bf Acknowledgment.} Part of this work is contained in the Doctoral thesis of the 
second author, who wish to thanks all the supporting institutions, especially the Doctoral School (ScuDo) of Politecnico of Torino. She wants especially thank Julianna 
Tymoczko for friendly and useful e-mail correspondence as well as hints by Prof. Knut- 
son, Prof. Lakshmibai and Prof. Raghavan. Both authors are deeply indebted with D. Laksov for substantial remarks and constant encouragement.

\section{Preliminaries and Notation}\label{preliminaries}

The presentation of the material below, basically covered in~\cite{tesi},  is strongly inspired by~\cite{LakTh} and~\cite{LakTh1}, where we learned to use  the polynomial $\ttp$ below.
\claim{}\label{recalling01}  Let  $A=\bigoplus_{i\geq 0}A_i$ be a graded ring such that $A_0=\ZZ$, $X$ an indeterminate over $A$, $M:=XA[X]$ and  $M(\ttp):=M/\ttp M$, where $\ttp=X^n+c_1X^{n-1}+\ldots+c_n\in A[X]$ is a monic polynomial of degree $n\geq 1$ such that $c_i\in A_i$.  For each $i\geq 1$, let $\ep^i:=X^i+\ttp M$. Then $M(\ttp)$ is a free $A$-module generated by  ${\bm\ep}=(\ep^1,\ldots,\ep^n)$. 
\claim{}\label{notpih} As in~\cite{anderson} (\S 5, Section~1) and~\cite{Gat1}, we prefer to use finite increasing sequences of integers instead of partitions. We shall indicate by ${\cal I}^k$ the set of all $I=(i_1,\ldots,i_k)\in\NN^k$ such that $1\leq i_1<\ldots<i_k$. The {\em weight} of $I\in{\cal I}^k$ is $wt(I)=\sum_{j=1}^k(i_j-j)$ (it is the weight of the corresponding partition $(i_k-k, i_{k-1}-(k-1),\ldots,i_1-1))$. 
If $I\in{\cal I}^k$, let
$
 {\cal P}(I):=\{H=(h_1,\ldots,h_k)
\in\NN^k\,\,|\,\,i_1\leq i_1+h_1<i_2\leq \ldots\leq i_{k-1}+h_{k-1}<i_k\}
$
and ${\cal P}(I,h)=\{(h_1,\ldots,h_k)\in{\cal P}(I)\,\,|\,\, \sum_{i=1}^kh_i=h\}$. Finally, let  ${\cal I}^k_n=\{I\in{\cal I}^k\,\,|\,\, i_k\leq n\}$.
 A routine check shows that if $I=(i_1,\ldots, i_k)\in {\cal I}^k_n$ and $H\in{\cal P}(I,h)$ then $I+H:=(i_1+h_1,\ldots,i_k+h_k)\in {\cal I}^k$. 

\claim{}\label{not23} Let  $\wM(\ttp):=\bigoplus_{k\geq 0}\wkM(\ttp)$ be the exterior algebra of $M(\ttp)$. If $I:=(i_1,\ldots,i_k)\in{\cal I}^k$, let $\w^I{\bm\ep}$ denote $\ikformep$. Each exterior power $\wkM(\ttp)$ is a free $A$-module with basis  $\Bw^k{\bm\ep}:=\{\w^I\ep:I\in{\cal I}^k_n\}$. If $a\in A_h$, the {\em weight} of $a\cdot\w^I{\bm\ep}$ is,  by definition,  $a+wt(I)$. 
Set $(\bw^kM({\tt p}))_w=\bigoplus_{0\leq h\leq w}\big(\bigoplus_{wt(I)=h}A_{w-h}\cdot\w^I{\bm \ep}\big)
$.
Then  $\wkM(\ttp)=\bigoplus_{w\geq 0}(\bw^k M({\tt p}))_w$,  a graded $A$-module via {\em weight}.
\claim{}  \label{rec24} By~\cite{GatSant1} there is a unique sequence $D:=(D_0,D_1,\ldots)$ (the {\em canonical ${\cal S}$-derivation}) of $A$-endomorphisms of $\wM(\ttp)$ such that i) the $h^{th}$-order {\em Leibniz's rule} 
\[
{ D_h(\alpha\w\beta)=\sum_{h_1,h_2\in\NN\,|\,h_1+h_2=h}D_{h_1}\alpha\w D_{h_2}\beta};
\]
holds for each $h\geq 0$ and each $\alpha,\beta\in\wkM(\ttp)$ and ii) satisfying the {\em initial conditions}  $D_h\ep^i=\ep^{i+h}$,  for each $h\geq 0$ and each   $ i\geq 1$.
 
\claim{} If $A[\TT]$ is the ring of polynomials in infinitely many indeterminates \linebreak $\TT:=(T_1,T_2,\ldots)$,  there is a natural evaluation map,  $\ev_D:A[\TT]\rightarrow End_A(\wM(\ttp))$,  sending $P\in A[\TT]$ to $P(D)$ (got by ``substituting"  $T_i\mapsto D_i$ into $P$). We denote by  ${\cal A}^*(\wM(\ttp))$ the image  of $\ev_D$ in $End_A(\wM(\ttp))$ and by ${\cal A}^*(\wkM(\ttp))$ the image of the natural restriction map  $\rho_k:{\cal A}^*(\wM(\ttp))\rightarrow End_A(\wkM(\ttp))$, sending  $P(D)$ to $P(D)_{|_{\wkM(\ttp)}}$. The following is a key result:
\begin{claim}{\bf Theorem.}\label{thm1rec} {\em The natural evaluation map $\ev_\kformep:{\cal A}^*(\wM(\ttp))\rightarrow\wkM(\ttp)$, mapping $P(D)\mapsto P(D)\kformep$ is surjective.}
\end{claim}

\proof The general determinantal formula by Laksov and Thorup (\cite{LakTh}) implies the following {\em Giambelli's formula}:
\be
\ikformep=\Delta_I(D)\cdot \kformep\label{eq:giambform},
\ee
where  $\Delta_I(D):=\ev_D(\Delta_I(\TT))\in {\cal A}^*(\wM(\ttp))$  and,  for each $I=(i_1,\ldots, i_k)\in\NN^k$, $\Delta_I(\TT):=\Delta_{(i_1,\ldots,i_k)}(\TT)=\det(T_{i_j-i})\in A[\TT]$ is the usual Schur polynomial in $(T_1,T_2,\ldots)$ (setting $T_0=1$ and $T_j=0$, if $j<0$).  \qed


\claim{} Theorem~\ref{thm1rec} easily implies that $\ker(\rho_k)=\ker\ev_\kformep$. The {\em Poincar\'e isomorphism}  $\Pi_k:{\cal A}^*(\wkM(\ttp))\rightarrow \wkM(\ttp)$ is  defined by $\rho_k(P(D))
\mapsto P(D)\kformep$.  In particular ${\cal A}^*(\wkM(\ttp))$ is a free $A$-module generated by $\{\Pi_k^{-1}(\w^I{\bm\ep})=\Delta_I(D):I\in{\cal I}^k_n\}$,  
 and $\wkM(\ttp)$ is a free ${\cal A}^*(\wkM(\ttp))$-module of rank $1$ generated by $\kformep$. 
 
\claim{} For each $k\geq 1$, let $\TT_k=(T_1,\ldots, T_k)$, so that $A[\TT_k]\subset A[\TT]$. Let $\widetilde{D}_i(\TT_k)$ ($i\geq 0$) be defined as:
$
\sum_{i\geq 0}\widetilde{D}_i({\bf T}_k)t^i:=\big(1+\Delta_{(2)}( \TT)t+\ldots+\Delta_{(2,\ldots, k+1)}(\TT)t^k\big)^{-1}\in A[[t]], 
$
and set
$
\widetilde{D}_i({\bf D}_k)=\ev_D(\widetilde{D}_i(\TT_k))=\widetilde{D}_i(D_1,\ldots, D_k)
$.
 One  has:
\bclm{\bf Theorem.}\label{generalpres} {\em The following  presentation holds:
\be
{\cal A}^*(\wkM({\tt p}))={A[D_1,\ldots, D_k]\over (\widetilde{D}_{n-k+1}({\bf D}_k,{\tt p}),\ldots, \widetilde{D}_{n}({\bf D}_k, \ttp))}.\label{eq:genpreshld}
\ee
where, for each $1\leq j\leq k$, 
$
\widetilde{D}_{n-k+j}({\bf D}_k,{\ttp})=\widetilde{D}_{n-k+j}({\bf D}_k)+\sum_{i=1}^{n-k+j}c_i\widetilde{D}_{n-k+j-i}({\bf D}_k),
$
}
\proof \cite{GatSant1}, Theorem~4.8.\qed
\eclm

\claim{} Due to the skew-symmetry of the exterior product, some cancellations occur  in the expansion of $D_h(\w^I{\bm\ep})\in \wkM(\ttp)$. The surviving summands are prescribed by Pieri's formula:
\be
D_h(\ep^{i_1}\w\ldots\w \ep^{i_k})=\sum_{{(H)\in {\cal P}(I,h)}}
 \ep^{i_1+h_1}\w\ldots \w \ep^{i_k+h_k}\label{eq:pieruno}
\ee
for each $\ikformep\in\wkM(\ttp)$ ($1\leq i_1<i_2<\ldots<i_k\leq n$). See~\cite{LakTh} or, since formula~(\ref{eq:pieruno}) is defined over the integers, use the same proof as in~\cite{Gat1}, Theorem~2.4.

\claim{} \label{prel211} Combining Pieri's formula~(\ref{eq:pieruno})  with Giambelli's-formula~(\ref{eq:giambform}), one has, for each $I\in{\cal I}^k$ and each $h\geq 0$:
\[
D_h\Delta_{I}(D)\kformep=D_h\cdot \w^I{\bm\ep}= \sum_{{H\in {\cal P}(I,h)}}
\w^{I+H}{\bm\ep}=\sum_{{H\in {\cal P}(I,h)}}
\Delta_{I+H}(D)\kformep,
\]
proving the equality  
$
D_h\Delta_{I}(D)=\sum_{{H\in {\cal P}(I,h)}}
\Delta_{I+H}(D)
$
in the ring  ${\cal A}^*(\wkM(\ttp))$.

\claim{}\label{fibcond} Let $A^*({\cal Y})$   and $H^*({\cal Y}):=H^*({\cal Y},\ZZ)$ be, respectively, the Chow intersection  and the integral cohomology ring of  a smooth complex projective variety ${\cal Y}$. Let \linebreak $p_k:G(k,E)\rightarrow {\cal Y}$ be the Grassmann bundle of $k$-planes of a rank $n$ vector bundle $p:E\rightarrow{\cal Y}$. Assume  that ${\cal Y}$ has a cellular decomposition (as in~\cite{Fu1}, p.~23).  Then, the natural cycle maps $A^*({\cal Y})\rightarrow H^*({\cal Y})$ and $A^*(G(k,E))\rightarrow H^*(G(k,E))$, doubling degrees, are ring isomorphisms (\cite{Fu1}, p.~378). 

\claim{} Let  $p:E\rightarrow {\cal Y}$  be a rank $n$-vector bundle as in~\ref{fibcond}. Let $A:=H^*({\cal Y})$,\linebreak 
$\ttp:=X^n+c_1X^{n-1}+\ldots+c_n$, where $c_i\in A_i:=H^{2i}({\cal Y})$ is the $i^{th}$ Chern class of $E$. The $A$-module structure of  $H^*(G(k,E))$ is defined by $a\cdot \alpha=p_k^*a\cup \alpha$, for each $a\in A$ and $\alpha\in H^*(G(k,E))$, where $\cup$ is the product in $H^*(G(k,E))$.  Let $0\rightarrow {\cal S}_k\rightarrow p_k^*E\rightarrow {\cal Q}_k\rightarrow 0$ be the tautological exact sequence over $G(k,E)$.  If $\xi=c_1(O_{\PP(E)}(1))$, then $\ttp(\xi)=0$,  by definition of Chern classes of $E$. By~\cite{Fu1}, p.~268, the ring $H^*(G(k,E))$ is a free  $A$-module with basis  $\{\Delta_I(c({\cal Q}_k-p_k^*E))\,|\; I\in{\cal I}^k_n\}$. \bclm{{\bf Theorem} (\cite{LakTh1}){\bf .}} {\em The $A$-module isomorphism 
$
\iota_k: A^*(G(k,E))\rightarrow {\cal A}^*(\wkM(\ttp))
$
defined by $\Delta_I(c({\cal Q}_k-p_k^*E))\mapsto \Delta_I(D)$ is a ring isomorphism.}
\eclm 

\proof It is enough to check on products of the form $c_h({\cal Q}_k-p_k^*E)\cup \Delta_I(c({\cal Q}_k-p_k^*E)$.
\begin{eqnarray*}
\iota_k(c_h({\cal Q}_k-p_k^*E)\cup \Delta_I(c({\cal Q}_k-p_k^*E))&=&\iota_k\big(\sum_{H\in {\cal P}(I,h)}\Delta_{I+H}(c({\cal Q}_k-p_k^*E))\big)=\\
=\sum_{H\in {\cal P}(I,h)}\Delta_{I+H}(D)=D_h\Delta_I(D)&=&\iota_k(c_h({\cal Q}_k-p_k^*E)\cdot \iota_k(\Delta_I((c({\cal Q}_k-p_k^*E)),
\end{eqnarray*}
by~\cite{Fu1}, Proposition~14.6.1, and~\ref{prel211}. \qed

\section{Equivariant Schubert Calculus} \label{sec3.0}
\claim{\bf Notation.} Let $T$ be a $d$-dimensional torus ($T$ may be either $(\CC^*)^d$ or $(S^1)^d$). For each $m\in\NN_{\geq 1}\cup\{\infty\}$, let $E_mT\rightarrow B_mT$ denote a principal $T$-bundle such that  $\pi_{i}(E_mT)=0$, for each $1\leq i\leq 2m$.  If $m<\infty$, $E_mT$ may be taken in fact, as we shall do, as the product of $d$-copies of $\CC^{m+1}\setminus\{{\bf 0}\}$ if $T=(\CC^*)^d$ and  as the product of $d$-copies of $(S^1)^{2m+1}$ if $T=(S^1)^d$ (\cite{anderson},~\cite{EdiGra}). In both cases,  $B_mT$ is the product of $d$-copies of a complex $\PP^{m}$, and  then it admits a cellular decomposition, meeting the hypothesis of~\ref{fibcond}.  If $m=\infty$, one shall simply write $ET\rightarrow BT$ instead of $E_\infty T\rightarrow B_\infty T$. It is a {\em universal principal $T$-bundle}, classifying principal $T$-bundles (up to homotopy) \cite{Brion}. 
\claim{} \label{comdiag32} Consider a $T$-action  on $\CC^n$  as in~\ref{int01}, ($*$). 
Let ${\cal E}_m=E_mT\times_T\CC^n:=E_m\times \CC^n/T$. Then $p_m:{\cal E}_m\rightarrow B_mT$ is the {\em associated vector bundle} to $E_mT\rightarrow B_mT$,  with fibers homeomorphic to $\CC^n$.
For each $1\leq k\leq n$ and each $m\geq 1$, $p_{k,m}:G(k,{\cal E}_m)\rightarrow B_mT$ is precisely the  associated bundle to $E_mT\rightarrow B_mT$ with fibers homeomorphic to $G(k,n)$, i.e.  $G(k,{\cal E}_m):=E_mT\times_T G(k,n)$.   For $m=\infty$ we shall simply write $p_k:G(k,{\cal E})\rightarrow BT$.
Let $m\geq m_1\geq m_2\geq 1$ three integers. The natural inclusions $E_{m_2}T\hookrightarrow E_{m_1}T\hookrightarrow E_mT$ and $B_{m_2}T\hookrightarrow B_{m_1}T\hookrightarrow B_mT$ induce the following diagram:
\[
\matrix{G(k,{\cal E}_{m_2})&\stackrel{f_{m_2,m_1}}\lra&G(k,{\cal E}_{m_1})&\stackrel{f_{m_1,m}}\lra&G(k,{\cal E}_m)\cr\downarrow&{}&\downarrow&{}&\downarrow\cr B_{m_2}T&\stackrel{g_{m_2,m_1}}{\lra}&B_{m_1}T&\stackrel{g_{m_1,m}}\lra&B_mT}
\]
with cartesian squares (Cf.~\cite{anderson}, \S 3, Lemma~1.4), where $f_{m_1,m}\circ f_{m_2,m_1}=f_{m_2,m}$ and $g_{m_1,m}\circ g_{m_2,m_1}=g_{m_2,m}$ and induces itself an arrow reversed diagram in cohomology:
\be
\matrix{H^*(G(k,{\cal E}_{m}))&\stackrel{\phi_{m,m_1}}\lra&H^*(G(k,{\cal E}_{m_1}))&\stackrel{\phi_{m_1,m_2}}\lra&H^*(G(k,{\cal E}_{m_2}))
\cr\uparrow&{}&\uparrow&{}&\uparrow\cr 
A(m)&\stackrel{\psi_{m,m_1}}{\lra}&A(m_1)&\stackrel{\psi_{m_1,m_2}}{\lra}&A(m_2)}\label{eq:eqdiag1}
\ee
where $\phi=f^*$, $\psi=g^*$ and   $A(m):=H^*(B_mT)$. Then $\phi_{m_1,m_2}\circ\phi_{m,m_2}=\phi_{m,m_2}$ and $\psi_{m_1,m_2}\circ\psi_{m,m_1}=\psi_{m,m_2}$. 
\claim{} \label{sec3.3} If $m_1>m_2$, the $A(m_2)$-module $H^*(G(k,{\cal E}_{m_2}))$ is a $2m_2$-{\em approximation} of the $A$-module $H^*(G(k,{\cal E}_{m_1}))$ in the following sense: the diagram
\be
\matrix{A(m_1)\otimes H^*(G(k,{\cal E}_{m_1}))&\lra &H^*(G(k,{\cal E}_{m_1}))\cr\psi_{m_1,m_2}\otimes\phi_{m_1,m_2}\Big\downarrow\,\,\,\,\,&{}&\psi_{m_1,m_2}\Big\downarrow\cr
A(m_2)\otimes H^*(G(k,{\cal E}_{m_2}))&\lra &H^*(G(k,{\cal E}_{m_2}))}\label{eq:eqdiag2}
\ee
commutes and the maps
\[
\left\{\matrix{{\phi_{m_1,m_2}}_{|_{H^i_T(G(k,n)}}&:&H^i_T(G(k,n))&\lra& H^i(G(k,{\cal E}_{m_2}))\cr
{\psi_{m_1,m_2}}_{|_{A_i}}&:&A(m_1)_i&\lra&A(m_2)_i}\right.
\]
are $\ZZ$-module isomorphisms for each $1\leq i\leq 2m_2$ (\cite{Brion}). Here $A(m)_i:=H^{2i}(B_mT)$.
\claim{} We shall write $\phi_m$ and $\psi_m$ instead of $\psi_{\infty,m}$ and $\phi_{\infty,m}$. By definition, \linebreak $H_T^*(G(k,n)):=H^*(G(k,{\cal E}))$ is the $T$-equivariant cohomology of the grassmannian $G(k,n)$. The $T$-equivariant cohomology $H_T^*(pt)$ of a point is  $A:=A(\infty)=H^*(BT)$. Diagram~(\ref{eq:eqdiag1}) says that $H_T^*(G(k,n))$ is an $A$-module which, by diagram~(\ref{eq:eqdiag2}), is $2m$-approximated by the $A(m)$-module $H^*(G(k,{\cal E}_m))$.
\claim{} Let  $F_k:=O_{G(k,n)}\otimes\CC^n$. All the terms occurring in the universal exact sequence over $G(k,n)$, $0\sra { T}_k\sra F_k\sra Q_k\sra 0$ ($S_k$ is the universal subbundle and $Q_k$ is the universal quotient bundle) are $T$-equivariant bundles (\cite{anderson}, \S 5, Section 2). Hence they possess {\em equivariant Chern classes} (\cite{Brion}, \cite{EdiGra}). If $F$ is an equivariant vector bundle, let $c^T(F):=\sum_{i\geq 0}c_i^T(F)t^i$ be the $T$-equivariant Chern polynomial. Define $c^T_t(Q_k-F_k)$ as $c_t(Q_k)/c^T_t(F_k)$ (the ratio taken in the ring $H_T^*(G(k,n))[[t]])$ and let $c^T(Q_k-F_k):=(c_i^T(Q_k-F_k))_{i\geq 0}$ be the sequence of its coefficients.  Let ${\cal Q}_{k,m}:=E_mT\times_TQ_k$ and notice that $p_{k,m}^*{\cal E}_m=E_mT\times_T  F_k$. General properties of equivariant Chern classes say that they pull back to Chern classes of the approximating bundles (see e.g.~\cite{anderson}, \S 4, Section~3), namely
$
\psi_m(c_i^T({ Q}_k))=c_i({\cal Q}_{k,m})$ { and} $\psi_m(c_i^T({F}_k))=c_i(p_{k,m}^*{\cal E}_{k,m})
$,
for each $1\leq i\leq 2m$. Denote by $\cup_T$ the product structure of $H_T^*(G(k,n))$.
\bclm{\bf Proposition.}\label{mainlemma} {\em The ring $H_T^*(G(k,n))$ is a free $A$-module generated by the classes $\{\delta_I:=\Delta_I(c^T(Q_k-F_k))\,|\, I\in{\cal I}^k_n\}$; in particular $(c_1^T(Q_k-F_k),\ldots,c_k^T(Q_k-F_k))$ generate it as an $A$-algebra and Pieri's formula holds:
\be
c_h^T(Q_k-F_k)\cup_T\delta_I=\sum_{H\in{\cal P}(I,h)}\delta_{I+H}.
\ee
}
\eclm
\proof Any $\alpha\in H_T^*(G(k,n))$ is a sum of homogeneous elements and then, without loss of generality, we may  assume  that $\alpha\in H_T^{2i}(G(k,n))$. Let $m>i$ and consider the $2m$-approximation
$
\phi_m:H_T^*(G(k,n))\lra H^*(G(k,{\cal E}_m)).
$
For each $I\in{\cal I}^k_n$, let  $\delta_{I,m}:=\Delta_I(c({\cal Q}_k-p_{k,m}^*{\cal E}_m)$. 
It follows that $\alpha_m:=\phi_m(\alpha)$ is a unique $H^*(B_mT)$-linear combination of $\delta_{I,m}$, with $wt(I)\leq i$ (\cite{Fu1}, p.~268, Proposition 14.6.5).  Then:
\[
\alpha=\phi_m^{-1}(\alpha_m)=\phi_m^{-1}(\sum_{wt(I)\leq i}a_{I,m}\delta_{I,m})=\psi_{m}^{-1}(a_{I,m})\delta_I\]
with $\psi_{m}^{-1}(a_{I,m})\in A$. An easy check shows that if $m_1,m_2>i$, then $\psi_{m_1}^{-1}(a_{I,m_1})=\psi_{m_2}^{-1}(a_{I,m_2})$ (by~\ref{comdiag32} and \ref{sec3.3}).
 This proves that $\{\delta_I\}$ generate $H_T^*(G(k,n))$ as $A$-module. Moreover, let $
\sum_{I\in{\cal I}^k_n}a_I\delta_I=0$ 
be any linear dependence relation. Without loss of generality we may assume that all the summands have the same degree $i=\deg(a_I)+wt(I)$ ($I\in{\cal I}^k_n$, see Section~\ref{notpih}). Then, for $m>i$,
\[
0=\phi_m(\sum_{I\in{\cal I}^k_n}a_I\delta_I)=\sum_{I\in{\cal I}^k_n}\psi_m(a_I)\phi_m(\delta_I)=\psi_m(a_I)\delta_{I,m},
\]
whence $\psi_m(a_I)=0$, for each $I\in{\cal I}^k_n$, because  $\{\delta_{I,m}:I\in{\cal I}^k_n\}$ is an $A(m)$-basis of $H^*(G(k,{\cal E}_m))$. Then $a_I=0$, for each $I\in{\cal I}^k_n$, by \ref{sec3.3}.
Similarly, for each $m>h+wt(I)$: 
\begin{eqnarray*}
&{}& c_h^T(Q_k-F_k)\cup_T\Delta_I(c^T(Q_k-F_k))=\phi_m^{-1}(\phi_m(c_h^T(Q_k-F_k)\cup_T\delta_I))=
\\
&{}&=\phi_m^{-1}(c_h({\cal Q}_{k,m}-p_{k,m}^*{\cal E}_m)\cup \delta_{I,m})=\phi_{m}^{-1}(\sum_{H\in{\cal P}(I,h)}\delta_{I+H,m})=\sum_{H\in{\cal P}(I,h)}\delta_{I+H}.\hskip 18pt\qed
\end{eqnarray*}

\noindent
We can finally state the main result of this paper.
\bclm{\bf Theorem.}\label{2ndmint} {\em The $A$-module isomorphism
$
\iota_k:H_T^*(G(k,n))\sra {\cal A}^*(\wkM(\ttp))
$
defined by $\delta_I\mapsto \Delta_I(D)$, is an $A$-algebras isomorphism. 
}
\eclm
\proof
The map $\iota_k$ is trivially a module isomorphism. One is then left to check that $\iota_k$ is indeed a ring isomorphism, i.e.  that  for each $I,J\in{\cal I}^k_n$, $\iota_k(\delta_J\delta_I)=\iota_k(\delta_J)\iota_k(\delta_I)$. Since any $\delta_J$ is a (Schur) polynomial in $(c_i^T(Q_k-F_k))$, it is enough to prove the claim for $J=(1,2,\ldots,k-1,k+h)$. One has:
\begin{eqnarray*}
\iota_k(c_h^T(Q_k-F_k)\cup_T\delta_I)&=&\iota_k 
(\sum_{H\in P(I,h)} \delta_{I+H})=\sum_{H\in {\cal P}(I,h)} \Delta_{I+H}(D)=\\=D_h\Delta_{I}(D)&=&
\iota_k(c_h^T(Q_k-F_k))\iota_k(\delta_I),
\end{eqnarray*}
by~\ref{prel211}, the definition of $\iota_k$ and~\ref{mainlemma}.\qed

\section{An Example}\label{KTex}
Let now $T$ be an $n$-dimensional torus, so that  $A:=H_T^*(pt)=\ZZ[y_1,\ldots,y_n]$. In this final section we shall apply Theorem~\ref{2ndmint} to the situation that  in~\cite{KT} is studied through {\em puzzles}. Very explicit {\em equivariant Pieri's formulas} will be also supplied (Theorem~\ref{TeoremaPieriEquivariante}).
\claim{}\label{KTex1}
 Let $A^{\oplus n}:=\bigoplus_{i=1}^{n}H_T^*(pt)$, with  the $A$-algebra structure given by  componentwise multiplication of polynomials. Denote by ${\frak S}^i$ the $i^{th}$ component of  ${\frak S}\in A^{\oplus n}$. 
If $T$ acts diagonally on $\CC^n$, then $H_T^*(\PP^{n-1})$ can be thought of as the  $A$-subalgebra of $A^{\oplus n}$,  generated, as $A$-module, by the {\em classes} ${\bf 1}, {\frak S}_1,\ldots,{\frak S}_{n-1}$, where ${\bf 1}^j=1_A$ and for all $1\leq i\leq n-1$,
$
{\frak S}^j_i=(y_j-y_1)\cdot\ldots\cdot(y_j-y_{i})=\prod_{h=1}^{i}(y_j-y_h)\in A,
$
imposing, according to  the recipe of~\cite{KT}, p.~230, the GKM conditions
(\cite{GKM}) to the components of ${\frak S}_i$.
For $1\leq i\leq n$, set $Y_i=y_i-y_1$: in particular $Y_1=0$.
Let 
$
\ttp=\prod_{j=1}^n(X-Y_j)\in A[X]
$,
$\ttp_0=1$ and  $\ttp_i=\prod_{1\leq j\leq i} (X-Y_j)$ for $1\leq i\leq n$. Clearly $\ttp_i$ divides $\ttp$ for each $1\leq i\leq n$ and $\ttp_n=\ttp$. 
An easy verification shows that
$
{\frak S}_i=\ttp_{i}({\frak S}_1):=\prod_{h=1}^{i}({\frak S}_1-Y_h{\bf 1}).
$
It is sufficient  to check it
for each component of ${\frak S}_i$:
\[
(\ttp_{i}({\frak S}_1))^j=\prod^{i}_{h=1}({\frak S}_1^j-(y_h-y_1))=\prod^{i}_{h=1}(y_j-y_1-y_h+y_1)=\prod^{i}_{h=1}(y_j-y_h)={\frak S}_i^j.
\]
Indeed, the $\ttp_i$ are a special case (for $G(1,n)$) of the {\em factorial Schur functions} (used in~\cite{Mihalcea2} and~\cite{formalism}).
Moreover $\ttp({\frak S}_1)=0$: in fact, for each $1\leq j\leq n$, at least a zero factor occurs into the product $\prod_{h=1}^n(y_j-y_h)$. We have hence checked that $H_T^*(\PP^{n-1})=A[X]/(\ttp)$ and that ${\frak S}_1=X+\ttp$. If  $M$ and $M(\ttp)$ are as in Section~\ref{recalling01},  there is a canonical $A$-algebra isomorphism
$\iota: H_T^*(\PP^{n-1})\sra {\cal A}^*(M(\ttp))$ sending ${\frak S}_1\mapsto D_1$. For each $i\geq 1$, let
\be
\mu^i=X\cdot\ttp_{i-1}+\ttp M\in M(\ttp)\label{eq:basemui}
\ee
Then ${\bm\mu}=(\mu^1,\ldots,\mu^n)$ is an $A$-basis of $M(\ttp)$. If $D:=(D_0,D_1,D_2,\ldots)$ is the canonical ${\cal S}$-derivation on $\wM(\ttp)$ (Cf.~\ref{rec24})  one concludes that
\be
D_1\mu^{j}=X(X\ttp_{j-1})+\ttp M=(X-Y_j)Xp_{j-1}+Y_jXp_{j-1}+\ttp M=\mu^{j+1}+Y_j\mu^{j}.\label{eq:defd1mui}
\ee

 \claim{} \label{sec4.3} The isomorphism $\iota_k$ between $H_T^*(G(k,n))$ and ${\cal
A}^*(\bigwedge^{k}(M(\ttp)))$ of Theorem~\ref{2ndmint} says that to study the former we can  work on the latter. The ring ${\cal
A}^*(\bigwedge^{k}(M(\ttp)))$ is a free $A$-module generated by $G_{I}(D)=\Pi_k^{-1}(\w^I{\bm\mu})$ (Cf.~\ref{not23} changing ${\bm\ep}$ into $\bm\mu$), while $H_T^*(G(k,n))$ is generated by ${\frak S}_I=\iota_k^{-1}(G_I(D))=\jmath_k^{-1}(\w^I{\bm\mu})$ (Cf.~\ref{seint02} for the notation).
A presentation of ${\cal A}^*(\wkM(\ttp))\cong H_T^*(G(k,n))$, in terms of generators and relations, is given by~(\ref{eq:genpreshld}), with $\ttp$ as in~\ref{KTex1} (see also~\cite{GatSant1}). Moreover, 
 since  $D_1,\ldots, D_k\in{\cal A}^*(\wM(\ttp))$  generate ${\cal A}^*(\wkM(\ttp))\cong H_T^*(G(k,n))$ as a ring (Cf.~(\ref{eq:genpreshld})),  an equivariant version of classical Pieri's formulas can be gotten by expressing  $D_l(\w^I{\bm\mu})$, for each $l\geq 0$ and each $I\in{\cal I}^k_n$, as an explicit $A$-linear combination of elements of $\bw^k{\bm\mu}$. To do that, we first need the following:


\begin{claim}{\bf Lemma.} {\em For all $i> 0$:
\begin{eqnarray} D_i(\mu^j)&=&D_1^i(\mu^j) = \sum_{l=0}^i
h_{i-l}(Y_{j},\ldots, Y_{j+l})\mu^{j+l},\quad  1\leq j\leq n  \label{eq:Dl}
\end{eqnarray}
where $h_0=1$ and, for each $m\geq 1$, 
 $
h_{m}(X_1,\ldots, X_k)$
is the complete homogeneous symmetric
polynomial in $X_1,\ldots,X_k$ of degree $m$ {\em (\cite{MacDon})}.}
\end{claim}
\proof The proof is by induction on the integer $i$. If $i=1$, one has, by~(\ref{eq:defd1mui})
\be
D_1\mu^j=\mu^{j+1}+Y_j\mu^j=
     \mu^{j+1}+h_1(Y_j)\mu^j. \label{eq:initstep}
\ee
Suppose that the formula holds for $i\geq 1$. Since
$D_{i}\mu^j=D_1D_1^{i-1}\mu^j$, by the inductive hypothesis:
\begin{eqnarray}
D_{i}(\mu^j)=D_1\bigg(\sum_{l=0}^{i-1}h_{i-l-1}(Y_{j},\ldots,
Y_{j+l})\mu^{j+l}\bigg)
=\sum_{l=0}^{i-1}h_{i-l-1}(Y_{j},\ldots, Y_{j+l})D_1\mu^{j+l}.\label{eq:initstep1}
\end{eqnarray}
\noindent Using~(\ref{eq:initstep}), last member of~(\ref{eq:initstep1}) is thence equal to:
\begin{eqnarray*}
&{}&\sum_{l=0}^{i-1}h_{i-l-1}(Y_{j},\ldots,
Y_{j+l})(\mu^{j+1}+h_1(Y_j)\mu^j)=\\
&=&\sum_{l=0}^{i-1}h_{i-l-1}(Y_{j},\ldots,Y_{j+l})\mu^{j+l+1}+
\sum_{l=0}^{i-1}h_{i-l-1}(Y_{j},\ldots,Y_{j+l})h_1(Y_{j+l})\mu^{j+l},
\end{eqnarray*}
\noindent
which may be rewritten as:
\begin{eqnarray*}
&=&\sum_{l'=1}^{i}h_{i-l'}(Y_{j},\ldots,Y_{j+l'-1})\mu^{j+l'}+\sum_{l=0}^{i-1}h_{i-l-1}(Y_{j},\ldots,Y_{j+l})h_1(Y_{j+l})\mu^{j+l}=\\
&=&h_{0}(Y_{j},\ldots,Y_{j+i-1})\mu^{j+i}+\sum_{l'=1}^{i-1}h_{i-l'}(Y_{j},\ldots,Y_{j+l'-1})\mu^{j+l'}+\\
&{}&+h_{i-1}(Y_{j})h_1(Y_{j})\mu^{j}+\sum_{l=1}^{i-1}h_{i-l-1}(Y_{j},\ldots,Y_{j+l})h_1(Y_{j+l})\mu^{j+l}.
\end{eqnarray*}
\noindent A further manipulation gives
\begin{eqnarray*}
&=&h_{0}(Y_{j},\ldots,Y_{j+i-1})\mu^{j+i}+h_{i}(Y_{j})\mu^{j}+ \\
&+&\sum_{l'=1}^{i-1}\bigg(h_{i-l'}(Y_{j},\ldots,Y_{j+l'-1})+h_{i-l'-1}(Y_{j},\ldots,Y_{j+l'})h_1(Y_{j+l'})\bigg)\mu^{j+l'}=\\
&=&h_{0}(Y_{j},\ldots,Y_{j+i-1})\mu^{j+i}+h_{i}(Y_{j})\mu^{j}+ \sum_{l'=1}^{i-1}h_{i-l'}(Y_{j},\ldots,Y_{j+l'})\mu^{j+l'}, 
\end{eqnarray*}
\noindent i.e.
$
D_{i}(\mu^j)= \sum_{l=0}^{i
}h_{i-l}(Y_{j},\ldots,Y_{j+l})\mu^{j+l}, \ \ 1\leq j\leq n,
$
as desired.\qed

\claim{} Let $I\in{\cal I}^k_n$, $\w^I{\bm\mu}\in\wkM(\ttp)$ and $l\geq 0$. Leibniz's rule  for $D_l$
gives:\linebreak
 $
D_l(\w^I{\bm\mu})=
\sum_{l_1+\ldots+l_k=l}D_{l_1}\mu^{l_1}\wedge\ldots\wedge
D_{l_k}\mu^{l_k}
$
(just by definition of $D_l$; see also~\cite{GatSant1}, Proposition~5.1, or~\cite{LakTh}).
 Using equation (\ref{eq:Dl}):
\begin{eqnarray*}
&{}&\sum_{l_1+\ldots+l_k=l}D_{l_1}\mu^{l_1}\wedge\ldots\wedge
D_{l_k}\mu^{l_k}=\\
&=&\sum_{l_1+\ldots+l_k=l}\left[\left(\sum_{m_{1}=0}^{l_1}h_{l_1-m_1}(Y_{i_1},Y_{i_1+1},\ldots,
Y_{i_1+m_1})\mu^{i_1+m_1}\right)\wedge\ldots\right.
\\
&& \qquad \qquad \qquad \qquad
\ldots\wedge\left.\left(\sum_{m_k=0}^{l_k}h_{l_{k}-m_k}(Y_{i_k},Y_{i_k+1},\ldots,
Y_{i_k+m_{k}})\mu^{i_k+m_{k}} \right)\right].
\end{eqnarray*}
Expanding the wedge products, last member can be written as:
\begin{eqnarray*}
&=&\sum_{l_1+\ldots+l_k=l}\left[\sum_{m_1+\ldots+m_k=0}^{l_1+\ldots+l_k}\left(\prod_{j=1}^k
h_{l_j-m_j}(Y_{i_j},Y_{i_j+1},\ldots, Y_{i_j+m_j})\mu^{i_1+m_1}
\wedge\ldots \wedge \mu^{i_k+m_{k}}\right)  \right]= \\
&=&\sum_{m_1+\ldots+m_k=0}^{l}\left[\left(\sum_{l_1+\ldots+l_k=l}\prod_{j=1}^k
h_{l_j-m_j}(Y_{i_j},Y_{i_j+1},\ldots, Y_{i_j+m_j})\right)
\mu^{i_1+m_1}
\wedge\ldots \wedge \mu^{i_k+m_{k}} \right]= \\
&=&\sum_{m_1+\ldots+m_k=0}^{l}\left[h_{l-\sum_{j=1}^km_j}(Y_{i_1},\ldots,
Y_{i_1+m_1},\ldots, Y_{i_k},\ldots, Y_{i_k+m_k}) \mu^{i_1+m_1}
\wedge\ldots \wedge \mu^{i_k+m_{k}} \right].
\end{eqnarray*}
Last equation   uses basic known properties of the complete symmetric polynomials (see e.g.~\cite{MacDon}). 
Putting $u=l-\sum_{j=1}^km_j$, one finally has:
\begin{center}
$
D_l(\mu^{i_1}\wedge\ldots\wedge\mu^{i_k})=
$
\end{center}
\be
=\sum_{u=0}^{l}\sum_{m_1+\ldots+m_k+u=l} h_{u}(Y_{i_1},\ldots,
Y_{i_1+m_1},\ldots, Y_{i_k},\ldots, Y_{i_k+m_k}) \mu^{i_1+m_1}
\wedge\ldots \wedge \mu^{i_k+m_{k}}.
\ee
The alternating feature of the $\wedge$-product causes cancellations of terms, so that, finally:

\begin{claim}\label{TeoremaPieriEquivariante} {\bf Theorem.} (See also~\cite{formalism}) {\em The following $T$-equivariant Pieri's formula holds:
\be
D_l(\w^I{\bm\mu})=\sum_{u=0}^{l}\sum_{M\in {\cal P}(I,l-u)} h_{u}(Y_{i_1},\ldots,
Y_{i_1+h_1},\ldots, Y_{i_k},\ldots, Y_{i_k+h_k})\cdot \w^{I+M}{\bm\mu},\label{eq:piereq}
\ee
where $M:=(m_1,\ldots,m_k)\in{\cal P}(I,l-u)$ is as in notation~\ref{notpih}.
}
\end{claim}

\proof (Cf.~\cite{Gat1}, Theorem~2.4).  By induction on the integer $k$. For $k=1$,
formula~(\ref{eq:piereq})  is trivially true.  First we prove it
directly for $k=2$. For each $l\geq 0$, let us split
sum~(\ref{eq:piereq}) as:
\begin{eqnarray*}
D_l(\mu^{i_1}\wedge\mu^{i_2})&=&\sum_{u=0}^{l}\sum_{m_1+m_2=l-u}
h_{u}(Y_{i_1},\ldots, Y_{i_1+m_1},Y_{i_1},\ldots,
Y_{i_2+m_2})\mu^{i_1+m_1} \wedge
\mu^{i_2+m_{2}}=\\
&=&\tt{U}+\overline{\tt{U}}
\end{eqnarray*}
 \[
{\rm where}\,\,\, {\tt U}=\sum_{u=0}^{l}\sum_{\scriptsize{\matrix{m_1+m_2=l-u\cr i_1+m_1<i_2}}}
h_{u}(Y_{i_1},\ldots, Y_{i_1+m_1},Y_{i_2},\ldots,
Y_{i_2+m_2})\mu^{i_1+m_1} \wedge \mu^{i_2+m_{2}}
 \]
\[
{\rm and}\,\,\, \overline{ {\tt U}}=\sum_{u=0}^{l}\sum_{\scriptsize{\matrix{m_1+m_2=l-u\cr
i_1+m_1\geq i_2}}} h_{u}(Y_{i_1},\ldots,
Y_{i_1+m_1},Y_{i_2},\ldots, Y_{i_2+m_2})\mu^{i_1+m_1} \wedge
\mu^{i_2+m_{2}}.
 \]
 One contends that $\overline{\tt U}$ vanishes. In fact,
on the finite set of all integers $i_2-i_1\leq a\leq l-u$,
define the bijection $\gamma(a)=i_2-i_1+l-u-a$. Then:
\begin{eqnarray*}
2\overline{\tt U}&=&\sum_{u=0}^{l}\sum_{m_1=i_2-i_1}^{l-u}
h_{u}(Y_{i_1},\ldots, Y_{i_1+m_1},Y_{i_2},\ldots,
Y_{i_2+l-u-m_{1}})\mu^{i_1+m_1} \wedge \mu^{i_2+l-u-m_{1}}+\\
&+&\sum_{u=0}^{l}\sum_{m_1=i_2-i_1}^{l-u} h_{u}(Y_{i_1},\ldots,
Y_{i_1+\gamma(m_1)},Y_{i_2},\ldots,
Y_{i_2+l-u-\gamma(m_1)})\mu^{i_1+\gamma(m_1)} \wedge
\mu^{i_2+l-u-\gamma(m_1)}=\\
&=&-\sum_{u=0}^{l}\sum_{m_1=i_2-i_1}^{l-u} h_{u}(Y_{i_1},\ldots,
Y_{i_1+m_1},Y_{i_2},\ldots,
Y_{i_2+l-u-m_{1}})\mu^{i_2+l-u-m_{1}} \wedge \mu^{i_2+m_1}+\\
&+&\sum_{u=0}^{l}\sum_{m_1=i_2-i_1}^{l-u} h_{u}(Y_{i_1},\ldots,
Y_{i_1+m_1},Y_{i_2},\ldots, Y_{i_2+l-u-m_1})\mu^{i_2+l-u-m_1}
\wedge \mu^{i_1+m_1}=0,
\end{eqnarray*}

\noindent hence $\overline{\tt U}=0$ and~(\ref{eq:piereq}) holds
for $k=2$. 

Suppose now that~(\ref{eq:piereq}) holds for all $1\leq
k'\leq k-1$. Then, for each $l\geq 0$:
\be
D_l(\mu^{i_1}\w\mu^{i_2}\w\ldots\w\mu^{i_k})=\sum_{{l'_k}+l_k=l}D_{l'_k}(\mu^{i_1}\wedge\ldots\wedge\mu^{i_{k-1}})\wedge
D_{l_k}\mu^{i_k}.\label{eq:nnpp}
\ee
\noindent By the inductive hypothesis, last member of~(\ref{eq:nnpp}) is equal to:
\be
\sum_{u=0}^{l'_k}{\cal H}_{k-1}(u)\w
\sum_{m_k=0}^{l_k}h_{l_{k}-m_k}(Y_{i_k},Y_{i_k+1},\ldots,
Y_{i_k+m_{k}})\mu^{i_k+m_{k}} ,\label{eq:piereq1}
\ee
where for notational brevity we set $I':=(i_1,\ldots,i_{k-1})$, $M'=(m_1,\ldots,m_{k-1})$ and:
\[
{\cal H}_{k-1}(u)=\sum_{M'\in{\cal P}(I',l'_k-u)} h_{u}(Y_{i_1},\ldots,
Y_{i_1+m_1},\ldots, Y_{i_{k-1}},\ldots, Y_{i_{k-1}+m_{k-1}})\w^{I'+M'}{\bm\mu},
\]

\noindent But now~(\ref{eq:piereq1})  can be equivalently written
as:
\be
\sum_{u=0}^{l-l''}{\cal H}_{k-2}(u)\wedge
D_{l''}(\mu^{i_{k-1}}\w\mu^{i_k})
\label{eq:pierieq2}
\ee
where
\[
{\cal H}_{k-2}(u)=\sum_{M''\in{\cal P}(I'',l-l"-u)}h_{u}(Y_{i_1},\ldots, Y_{i_1+m_1},\ldots,
Y_{i_{k-2}},\ldots, Y_{i_{k-2}+m_{k-2}}) \mu^{i_1+m_1}
\wedge\ldots \wedge \mu^{i_{k-2}+m_{k-2}},
\]
and where $M''=(m_1,\ldots, m_{k-2})$.
The inductive hypothesis implies that
\begin{small}
\begin{eqnarray*}
&{}&D_{l''}(\mu^{i_{k-1}}\wedge\mu^{i_k})=\\
&=&\sum_{u=0}^{l''}\sum_{\scriptsize{\matrix{m_{k-1}+m_k=l''-u\cr
i_{k-1}+m_{k-1}<i_k}}} h_{u}(Y_{i_{k-1}},\ldots,
Y_{i_{k-1}+m_{k-1}},Y_{i_k},\ldots,
Y_{i_k+m_k})\mu^{i_{k-1}+m_{k-1}} \wedge \mu^{i_k+m_{k}}
\end{eqnarray*}
\end{small}

 \noindent which, substituted
into~(\ref{eq:pierieq2}), yields precisely~(\ref{eq:piereq}).\qed

\bclm{\bf Example.} i) The coefficient of $\mu^2\w\mu^3\w\mu^7$, in
the expansion of $D_3(\mu^2\w\mu^3\w\mu^5)$, is
$h_1(Y_2,Y_3,Y_5,Y_6,Y_7)=Y_2+Y_3+Y_5+Y_6+Y_7=y_2+y_3+y_5+y_6+y_7-5y_1
$; ii) for 

\noindent
$l=1$,  Pieri's formula~(\ref{eq:piereq})  reads:
\be
D_1(\mu^{i_1}\wedge\ldots\wedge\mu^{i_k})=\sum_{j=1}^k\mu^{i_1}\wedge\ldots\wedge\mu^{i_j+1}\w \ldots\w
\mu^{i_k}+(Y_{i_1}+\ldots+Y_{i_k})\mu^{i_1}\wedge\ldots\w
\mu^{i_k}.\label{eq:piericod1}
\ee
\eclm
\claim{} \label{ex55} Let $1\leq k\leq n $ be  a fixed integer and
 $\left\{{\tiny \matrix{n \cr k}}\right\}$  the set of all functions \linebreak ${\bf a}_{i_1\ldots i_{k}}:\{1,\ldots,n\}\sra
\{0,1\}$ ($1\leq i_1<\ldots<i_k\leq n$) such that ${\bf a}_{i_1\ldots i_{k}}(j)=0$ if $j\in \{i_1,\ldots,i_{k}\}$ and ${\bf a}_{i_1\ldots i_{k}}(j)=1$  otherwise (see~\cite{KT}). Each $\lambda\in\left\{{\tiny \matrix{n \cr k}}\right\}$ can be also represented by a string $\lambda(1)\lambda(2)\ldots\lambda(n)$ of zeros and ones only.
If $\lambda={\bf a}_{i_1\ldots i_k}$, write $\SS_\lambda:=\Pi_k^{-1}(\ikformu)\in{\cal A}^*(\wkM(\ttp))$  and $\SS_{div}:=\SS_{{\bf a}_{12\ldots k-1,k+1}}=\SS_{0\ldots010}$. It is easily seen that $
\tilde{S}_{div}(D)=\left(D_1-\sum_{r=1}^k Y_r\right)$ and then, using~(\ref{eq:defd1mui}) and~(\ref{eq:piericod1}):
\begin{small}
\begin{eqnarray}
&&\SS_{div}\SS_\lambda\kformu=\SS_{div}\mu^{i_1}\wedge\ldots\wedge\mu^{i_k}=\left(D_1-\sum_{r=1}^k Y_r\right)\ikformu= \nonumber \\
&=&D_1\ikformu-\sum_{r=1}^k Y_r\ikformu=\nonumber\\
&=&\sum_{j=1}^k\mu^{i_1}\wedge\ldots\wedge\mu^{i_{j-1}}\wedge\mu^{i_j+1}\w\mu^{i_{j+1}}\w
\ldots\w
\mu^{i_k}+\left(\sum_{r=1}^k(Y_{i_r}-Y_r)\right)\mu^{i_1}\wedge\ldots\w
\mu^{i_k}=\nonumber\\
&=&\sum_{j=1}^k\mu^{i_1}\wedge\ldots\wedge\mu^{i_{j-1}}\wedge\mu^{i_j+1}\w\mu^{i_{j+1}}\w
\ldots\w
\mu^{i_k}+\left(\sum_{r=1}^k(y_{i_r}-y_r)\right)\mu^{i_1}\wedge\ldots\w
\mu^{i_k}.\nonumber\label{PieriTGKT}
\end{eqnarray}
\end{small}
\noindent
We have hence proven that, as in~\cite{KT}, 
$\SS_{div}\SS_\lambda=\sum_{\lambda':\lambda'\rightarrow
\lambda}\widetilde{S}_{\lambda'}+\left(\widetilde{S}_{div}|{\lambda}\right)\widetilde{S}_{\lambda}.
$
The notation  $\lambda':\lambda'\rightarrow
\lambda$ means that $\lambda'$ differs from $\lambda$ in only two
spots $i, i+1$, where $\lambda$ has $01$ and $\lambda'$ has $10$ (\cite{KT}, p.~236, Lemma 3),  while the coefficient $\left(\widetilde{S}_{div}|{\lambda}\right)$
is given by
$
\left(\widetilde{S}_{div}|{\lambda}\right)=
\sum_{j=1}^{n}(1-\lambda(j))y_j-\sum_{i=1}^{k}y_i
$. 
Thus, if  $\lambda={\bf a}_{i_1,\ldots,i_k}$, one has, in fact:
\begin{eqnarray*}
\left(\widetilde{S}_{div}|{\lambda}\right)&=&
(y_{i_1}+y_{i_2}+\ldots+y_{i_k})-(y_1+y_2+\ldots+y_k)=\sum_{r=1}^k
(y_{i_r}-y_{r}).
\end{eqnarray*}

\claim{\bf Remark.} 
In~\cite{formalism} and~\cite{formalism1}  Laksov  makes  explicit the exact relationship between the basis $\bw^k{\bm\mu}$ and that of Knutson and Tao. Moreover, using the deformed polynomial $\ttp(q)=\prod_{i=1}^n(X-Y_i)+q\in A[q]$, he recovers Mihalcea's work (\cite{Mihalcea1}, \cite{Mihalcea2})  within the framework of~\cite{LakTh} and~\cite{LakTh1}: the small (equivariant) quantum cohomology of $G(1,n)$ is well behaved with respect to ``taking exterior powers" (in spite of its non functoriality)!

\claim{\bf Example. }\label{KTrevisited} (Cf.~\cite{KT}, p.~231).   Let $A:=\ZZ[y_1,y_2,y_3,y_4]$, $Y_i=y_i-y_1$ ($1\leq i\leq 4$),
$
\ttp=X(X-Y_2)(X-Y_3)(X-Y_4)
$,
$M=XA[X]$ and $M(\ttp)=M/\ttp M$. Let ${\bm\mu}=(\mu^1,\mu^2,\mu^3,\mu^4)$ as in~(\ref{eq:basemui}).
Then
\[
\bw^2{\bm\mu}=(\mu^1\w\mu^2,\quad \mu^1\w\mu^3,\quad \mu^1\w\mu^4,\quad \mu^2\w\mu^3,\quad \mu^2\w\mu^4,\quad \mu^3\w\mu^4)
\]
is the basis of $\bw^2M(\ttp)$ (degree and lexicographically (totally) ordered) induced by ${\bm\mu}$. Let $G_{ij}(D)={\Pi}_2^{-1}(\mu^i\w\mu^j)\in{\cal A}^*(\Bw^2M(\ttp))$.  ``Integration by parts" (as in~\cite{GatSant1}) or the general  determinantal formula due to Laksov and Thorup (\cite{LakTh}), yields:
\[
G_{12}(D)=1;\qquad
G_{13}(D)=D_1-Y_2D_0;
\]
\[
G_{14}(D)=D_2-e_1(Y_2, Y_3)D_1+e_2(Y_2,Y_3)D_0;\qquad
 G_{23}(D)=D_1^2-D_2;
\]
\[
G_{24}(D)=(D_1-Y_4)G_{14}(D);\qquad
G_{34}(D)=(D_2-Y_4^2)G_{14}^{\bm
\mu}(D)-(Y_2+Y_4)G_{24}(D).
\]
Let us write the matrices associated to $G_{ij}(D)$ in the
basis $\wedge^2{\bm\mu}$. The matrix $(G_{12}(D))$ is just the $6\times 6$ identity matrix.
We show, for instance, how the computation of  $(G_{13}(D))$ goes on. Let us express $G_{13}(D)\cdot \w^I{\bm\mu}$ as an $A$-linear combination of $\bw^2{\bm\mu}$. One has, e.g.:
\begin{small}
\begin{eqnarray*}
\bigg(D_1-Y_2D_0\bigg)\mu^1\w \mu^2 &=& \mu^2\w \mu^2+\mu^1\w(\mu^3+Y_2\mu^2)-Y_2\mu^1\w\mu^2 = \mu^1\w \mu^3
\end{eqnarray*} 
and 
 \begin{eqnarray*}
  G_{13}(D)\mu^1\w\mu^3&=&\bigg(D_1-Y_2D_0\bigg)\mu^1\w \mu^3=\mu^2\w \mu^3+\mu^1\w(\mu^4+Y_3\mu^2)-Y_2\mu^1\w\mu^3 =\\
                                         &=& \mu^2\w \mu^3+ \mu^1\w\mu^4+{(Y_3-Y_2)}\mu^1\w\mu^3.                                      
\end{eqnarray*}                             
   \end{small}                                      
Continuing in the same way,  the matrix of $G_{13}(D)$ in the basis $\bw^2{\bm\mu}$ is:
\begin{scriptsize}
\begin{eqnarray*}
(G_{13}(D))=\left(\matrix{
 {\begin{tabular}{|l|}
\hline
{ $0$}\\
\hline
\end{tabular}}& 0    &0& 0  &  0  & 0\cr
  1 & {\begin{tabular}{|l|}
\hline
{ $Y_3-Y_2$}\\
\hline
\end{tabular}}
&0& 0  &  0  & 0\cr
  0 & 1 & {\begin{tabular}{|l|}
\hline
{ $Y_4-Y_2$}\\
\hline
\end{tabular}} & 0  &  0  & 0\cr
  0 &0&0&{\begin{tabular}{|l|}
\hline
{ $Y_3$}\\
\hline
\end{tabular}}& 0  & 0 \cr
  0 &0 & 1& 0  & {\begin{tabular}{|l|}
\hline
{ $Y_4$}\\
\hline
\end{tabular}} & 0\cr
  0 &0&0& 0&  1  &{\begin{tabular}{|l|}
\hline
{ $Y_4+Y_3-Y_2$}\\
\hline
\end{tabular}}
}%
\right)
\end{eqnarray*}
\end{scriptsize}
Similarly, one can compute (by hands) the matrices of the remaining $G_{ij}(D)$, getting
\begin{scriptsize}
\begin{eqnarray*}
(G_{14}(D))=\pmatrix{
  {\begin{tabular}{|l|}
\hline
{ $0$}\\
\hline
\end{tabular}}
&    0      &         0       &     0     &          0            & 0\cr
       0     &
{\begin{tabular}{|l|}
\hline
{ $0$}\\
\hline
\end{tabular}}
&         0       &     0     &          0            & 0\cr
       1     & Y_4-Y_2   & {\begin{tabular}{|l|}
\hline
{ $(Y_4-Y_2)(Y_4-Y_3)$}\\
\hline
\end{tabular}}&     0     &          0            & 0\cr
       0     &    0      &         0       &{\begin{tabular}{|l|}
\hline
{ $0$}\\
\hline
\end{tabular}}&          0            & 0 \cr
       0     &    1      &      Y_4-Y_3    &     Y_4   & {\begin{tabular}{|l|}
\hline
{ $Y_4(Y_4-Y_3)$}\\
\hline
\end{tabular}} & 0\cr
       0     &    0      &         1       &     0     &         Y_4           &{\begin{tabular}{|l|}
\hline
{ $Y_4(Y_4-Y_2)$}\\
\hline
\end{tabular}}}%
\end{eqnarray*}

\begin{eqnarray*}
(G_{23}(D))=\left(\matrix{
 {\begin{tabular}{|l|}
\hline
{ $0$}\\
\hline
\end{tabular}}
 &    0      &         0       &     0     &          0            & 0\cr
       0     &{\begin{tabular}{|l|}
\hline
{ $0$}\\
\hline
\end{tabular}}&         0       &     0     &          0            & 0\cr
       0     &    0      & {\begin{tabular}{|l|}
\hline
{ $0$}\\
\hline
\end{tabular}}    &     0     &          0            & 0\cr
       1     &   Y_3     &         0       &{\begin{tabular}{|l|}
\hline
{ $Y_2Y_3$}\\
\hline
\end{tabular}}
&          0            & 0 \cr
       0     &    1      &        Y_4      &    Y_2    &{\begin{tabular}{|l|}
\hline
{ $Y_2Y_4$}\\
\hline
\end{tabular}}
& 0\cr
       0     &    0      &         0       &     1     &         Y_4           &{\begin{tabular}{|l|}
\hline
{ $Y_3Y_4$}\\
\hline
\end{tabular}}}%
\right)
\end{eqnarray*}

\begin{eqnarray*}
(G_{24}(D))=\left(\matrix{
{\begin{tabular}{|l|}
\hline
{ $0$}\\
\hline
\end{tabular}}
 &    0      &         0 &     0     &          0            & 0\cr
       0     & {\begin{tabular}{|l|}
\hline
{ $0$}\\
\hline
\end{tabular}}
&         0 &     0     &          0            & 0\cr
       0     &    0      &{\begin{tabular}{|l|}
\hline
{ $0$}\\
\hline
\end{tabular}}
&     0     &          0            & 0\cr
       0     &    0      &         0 &{\begin{tabular}{|l|}
\hline
{ $0$}\\
\hline
\end{tabular}}
&          0            & 0 \cr
       1     &   Y_4     &Y_4(Y_4-Y_3)&    Y_2Y_4     & {\begin{tabular}{|l|}
\hline
{ $Y_2Y_4(Y_4-Y_3)$}\\
\hline
\end{tabular}}
& 0\cr
       0     &    1      &      Y_4   &     Y_4       &         Y_4^2           &{\begin{tabular}{|l|}
\hline
{ $Y_4Y_3(Y_4-Y_2)$}\\
\hline
\end{tabular}}
}%
\right)
\end{eqnarray*}

\begin{eqnarray*}
(G_{34}(D))=\pmatrix{
 {\begin{tabular}{|l|}
\hline
{ $0$}\\
\hline
\end{tabular}}
 &    0      &         0       &     0     &        0  & 0\cr
       0     & {\begin{tabular}{|l|}
\hline
{ $0$}\\
\hline
\end{tabular}}
&         0       &     0     &        0  & 0\cr
       0     &    0      &{\begin{tabular}{|l|}
\hline
{ $0$}\\
\hline
\end{tabular}}
&     0     &        0  & 0\cr
       0     &    0      &         0       &{\begin{tabular}{|l|}
\hline
{ $0$}\\
\hline
\end{tabular}}
&        0  & 0 \cr
       0     &    0      &         0       &     0     &{\begin{tabular}{|l|}
\hline
{ $0$}\\
\hline
\end{tabular}}
& 0\cr
       1     &Y_4+Y_3-Y_2& Y_4(Y_4-Y_2)    &   Y_4Y_3  &Y_4Y_3(Y_4-Y_2)&{\begin{tabular}{|l|}
\hline
{ $Y_3(Y_4-Y_2)Y_4(Y_3-Y_2)$}\\
\hline
\end{tabular}}
}
\end{eqnarray*}
\end{scriptsize}

\noindent  
 Let $\{\SS_{{\bf a}_{ij}}\,|\, (i,j)\in{\cal I}^2_4\}$ be the basis of $H_T^*(G(2,4))$ depicted in~\cite{KT}, p.~231, Fig.~7. Then the inverse ${\cal A}^*(\bw^2M(\ttp))\sra H_T^*(G(2,4))$ of the ring isomorphism $\iota_2$ (as  in~\ref{sec4.3}, for $k=2$) is explicitly given by $G_{ij}(D)\mapsto \SS_{{\bf a}_{ij}}$. Notice that the boxed diagonal entries of the matrices  $(G_{ij}(D))$ satisfy the GKM conditions (\cite{KT}). 

\medskip

Dipartimento di Matematica, 

Politecnico di Torino, 

C.so Duca degli Abruzzi, 24, 

10129, Torino

\medskip
{\tt letterio.gatto@polito.it}

{\tt taise@calvino.polito.it}

 \end{document}